%% file: ProjectiveQuiver.tex
\documentclass{amsart}

\usepackage{graphicx}
\usepackage{a4wide}
\usepackage{psfrag}
\usepackage[all]{xy}

\let\cal\mathcal
\def\AA{{\cal A}}

\def\CC{{\cal C}}

\def\II{{\cal I}}

\def\PP{{\cal P}}
\def\QQ{{\cal Q}}
\def\RR{{\cal R}}

\let\blb\mathbb

\def\bZ{{\blb Z}}

\def\bN{{\blb N}}

\def\bZ{{\blb Z}}

\def\rad{\operatorname {rad}}

\def\rep{\operatorname{rep}}

\def\Ext{\operatorname {Ext}}

\def\Hom{\operatorname {Hom}}

\def\Sl{\operatorname {Sl}}

\def\infrad{\operatorname {rad}^{\infty}}

\def\r{\rightarrow}
\def\l{\leftarrow}

\def\d{\downarrow}

\DeclareMathOperator{\ind}{ind}

\newcommand\Db{D^{b}}
\newcommand\Kb{K^{b}}

\renewcommand\r{d^{\bullet}}
\renewcommand\l{d_{\bullet}}

\newcommand\rQ{d^{\bullet}_{Q}}
\renewcommand\d{d}
\newcommand\dQ{d_Q}

\newcommand\Sr[2]{S^{\bullet}(#1,#2)}
\renewcommand\Sl[2]{S_{\bullet}(#1,#2)}

\newcommand\tm{\tau^{-1}}
\renewcommand\t{\tau}

\newtheorem{lemma}{Lemma}[section]
\newtheorem{proposition}[lemma]{Proposition}
\newtheorem{theorem}[lemma]{Theorem}
\newtheorem{corollary}[lemma]{Corollary}

\theoremstyle{definition}

\newtheorem{example}[lemma]{Example}

{

}

\theoremstyle{remark}

\newdimen\uboxsep \uboxsep=1ex
\def\uboxn#1{\vtop to 0pt{\hrule height 0pt depth 0pt\vskip\uboxsep
\hbox to 0pt{\hss #1\hss}\vss}}

\def\uboxs#1{\vbox to 0pt{\vss\hbox to 0pt{\hss #1\hss}
\vskip\uboxsep\hrule height 0pt depth 0pt}}

\newcommand\exa{\nopagebreak \begin{center}\smallskip \nopagebreak               \begin{minipage}[t]{6cm}\sloppy}
\newcommand\exb{\end{minipage}\kern 1cm\begin{minipage}[t]{8cm}\sloppy}
\newcommand\exc{\end{minipage}\kern -3cm \smallskip\end{center}}

\title{The Quiver of Projectives in Hereditary Categories with Serre duality}
\author{Carl Fredrik Berg}
\author{Adam-Christiaan van Roosmalen}
\address{Carl Fredrik Berg\\Institutt for matematiske fag\\
NTNU\\7491 Trondheim\\Norway\\
(Currently working for StatoilHydro R\&D Centre\\Arkitekt Ebbells veg 10\\Rotvoll\\7053 Ranheim\\Norway)} \email{carlpaatur@hotmail.com}
\address{Adam-Christiaan van Roosmalen\\Max-Planck-Institut f\"{u}r Mathematik
\\Vivatsgasse 7\\53111 Bonn\\Germany}\email{vroosmal@mpim-bonn.mpg.de}

\begin{document}

\begin{abstract}
Let $k$ be an algebraically closed field and $\AA$ a $k$-linear hereditary category satisfying Serre duality with no infinite radicals between the preprojective objects.  If $\AA$ is generated by the preprojective objects, then we show that $\AA$ is derived equivalent to $\rep_k Q$ for a so called \emph{strongly locally finite} quiver $Q$.  To this end, we introduce light cone distances and round trip distances on quivers which will be used to investigate sections in stable translation quivers of the form $\bZ Q$.
\end{abstract}

\maketitle

\tableofcontents

\input{Introduction}
\input{Preliminaries}
\input{Distances}
\input{Main}
\input{Categories}

\bibliographystyle{amsplain}

\providecommand{\bysame}{\leavevmode\hbox to3em{\hrulefill}\thinspace}
\providecommand{\MR}{\relax\ifhmode\unskip\space\fi MR }
\providecommand{\MRhref}[2]{%
  \href{http://www.ams.org/mathscinet-getitem?mr=#1}{#2}
}
\providecommand{\href}[2]{#2}

\end{document}

%% file: Introduction.tex
\section{Introduction}

With a quiver $Q$ we may associate a stable translation quiver $\bZ Q$ as follows: the vertices are given by $(n,x)$ where $n \in \bZ$ and $x \in Q$.  The number of arrows $(i,x) \to (j,y)$ is equal to the number of arrows $x \to y$ in $Q$ if $i=j$, equal to the number of arrows $y \to x$ if $j = i+1$, and equal to zero otherwise.  On the vertices of $\bZ Q$, we may define a translation $\t : \bZ Q \to \bZ Q$ by $\t(n,x) = (n-1,x)$.  This is an automorphism of $\bZ Q$ that makes $\bZ Q$ a stable translation quiver.

Non-isomorphic quivers $Q$ and $Q'$ may give rise to isomorphic stable translation quivers $\bZ Q$ and $\bZ Q'$.  We define a section of $\bZ Q$ as a full subquiver $Q'$ of $\bZ Q$ such that the embedding $Q' \to \bZ Q$ extends to an isomorphism $\bZ Q' \to \bZ Q$ of stable translation quivers.  In this paper, we will investigate for which quivers $Q$ the stable translation quiver $\bZ Q$ admits a \emph{strongly locally finite} section $Q'$, i.e. every vertex of $Q'$ has finitely many neighbors and $Q'$ is without subquivers of the form $\cdot \to \cdot \to \cdots$ or $\cdots \to \cdot \to \cdot$

Before stating our main result, we will need a definition.  Let $Q$ be a quiver.  For two vertices $x,y \in Q$ we define the \emph{round trip distance} $\d(x,y)$ as the least number of arrows that have to be traversed in the opposite direction on an unoriented path from $x$ to $y$ and back to $x$.  If $Q$ does not have oriented cycles, then for all $x,y,z \in Q$
\begin{enumerate}
  \item $\d(x,y) \geq 0$ and $\d(x,y)=0 \Longleftrightarrow x=y$,
  \item $\d(x,y) = \d(y,x)$,
  \item $\d(x,z) \leq \d(x,y) + \d(y,z)$
\end{enumerate}
such that $\d$ defines a distance on the vertices of $Q$ (Proposition \ref{proposition:d_is_distance} in the text).  To the round trip distance, we may associate \emph{round trip distance spheres} as follows
$$S(x,n) = \{y \in Q \mid \d(x,y) = n\}.$$

We may now formulate our main theorem (Theorem \ref{theorem:Main} in the text).

\begin{theorem}\label{theorem:Introduction}
Let $Q$ be a connected quiver, then the following are equivalent.
\begin{itemize}
\item
The quiver $Q$ has no oriented cycles, and for a certain $x \in Q$ (or equivalently: for all $x \in Q$) the round trip distance spheres $S_Q(x,n)$ are finite, for all $n \in \mathbb{N}$.
\item
There are only finitely many paths in $\bZ Q$ between two vertices.
\item
The translation quiver $\mathbb{Z} Q$ has a strongly locally finite section.
\end{itemize}
\end{theorem}

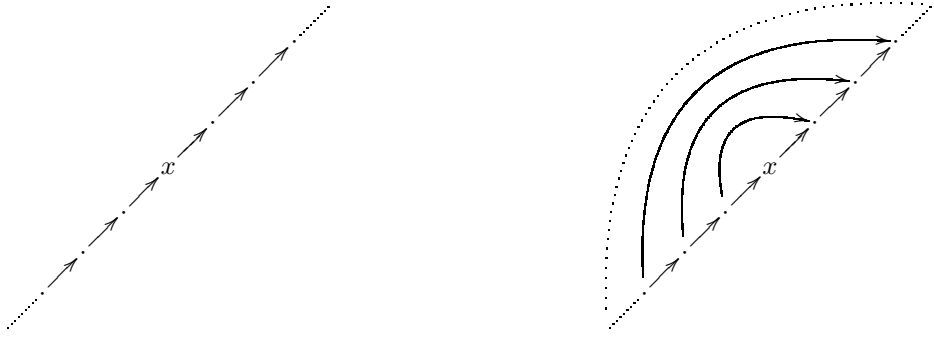
\begin{figure}[tb]
	\centering
		\exa 
		$$\xymatrix@1@C=15pt@ur{\ar@{.}[r]&\cdot\ar[r]&\cdot\ar[r]&\cdot\ar[r]&x\ar[r]&\cdot\ar[r]&\cdot\ar[r]&\cdot\ar@{.}[r]&}$$
		\exb
		$$\xymatrix@1@C=15pt@R=10pt@ur{ \ar@{.}[r] \ar@{.}@/^50pt/[rrrrrrrr] & \cdot \ar[r]
\ar@/^40pt/[rrrrrr] & \cdot \ar[r] \ar@/^30pt/[rrrr] & \cdot \ar[r]
\ar@/^20pt/[rr] & x \ar[r] & \cdot \ar[r] & \cdot \ar[r] & \cdot
\ar@{.}[r] & }$$
		\exc
	\caption{A quiver satisfying the equivalent conditions of Theorem \ref{theorem:Introduction} (left) and one that does not (right)}
	\label{figure:ExampleQuivers}
\end{figure}

As an example, we see that the left quiver in Figure \ref{figure:ExampleQuivers} satisfies the first condition of the previous theorem, while in the quiver on the right hand side the round trip distance sphere $S_Q(x,1)$ has infinitely many vertices.

Our main reason to investigate this problem has been a question by Reiten and Van den Bergh in \cite{ReVdB02}.  In that article, Reiten and Van den Bergh classified all $k$-linear noetherian abelian hereditary Ext-finite categories with Serre duality.  One type of such categories, characterized by being generated by preprojectives, was constructed by formally inverting a right Serre functor in the category $\rep_k Q$ of finitely presented representations of a certain quiver $Q$.

Reiten and Van den Bergh suggest another construction of these categories, and a shorter proof of their classification, based on the answer to the following question.  Let $\AA$ be an hereditary noetherian category with Serre duality, and let $Q$ be the full subquiver of the Auslander-Reiten quiver of $\AA$ spanned by the isomorphism classes of the indecomposable projectives.  Does $\bZ Q$ have a strongly locally finite section?

Since Reiten and Van den Bergh note (\cite[Lemma II.3.1]{ReVdB02}) that for the quivers $Q$ under consideration there are only finitely many paths between two vertices in $\bZ Q$, our Theorem \ref{theorem:Introduction} gives a positive answer.  Following the ideas of Reiten and Van den Bergh, we obtain an alternative way of constructing the noetherian categories generated by preprojectives (Ringel already gave an alternative way of constructing such categories using ray quivers in \cite{Ringel02b}).

\begin{theorem}
Let $\AA$ be a noetherian $k$-linear abelian Ext-finite hereditary category with Serre duality.  Assume $\AA$ is generated by the preprojective objects, then $\AA$ is derived equivalent to $\rep Q'$ where $Q'$ is strongly locally finite.
\end{theorem}

A slightly more general result, not involving the noetherian condition, is given by Theorem \ref{theorem:MainDerived}.

Let $\AA$ be a hereditary category with Serre duality.  The following theorem (Corollary \ref{corollary:Classification} in the text) characterizes all quivers which occur as a subquiver of the Auslander-Reiten quiver of $\AA$ generated by indecomposable projectives (called the \emph{quiver of projectives}).  This complements a result from \cite{ReVdB02} where all such quivers that can arise when $\AA$ is noetherian were characterized as being star quivers.

\begin{theorem}
Let $Q$ be a quiver.  There is an abelian hereditary category with Serre duality having $Q$ as its quiver of projectives if and only if $Q$ satisfies the equivalent conditions of Theorem \ref{theorem:Introduction}.
\end{theorem}

The proof of Theorem \ref{theorem:Introduction} is a constructive one.  Let $Q$ be a quiver.  In $\bZ Q$ we define the \emph{right light cone} centered on a vertex $x \in \bZ Q$ as the set of all vertices $y$ such that there is an oriented path from $x$ to $y$ but not to $\t y$.  Dually, we define the \emph{left light cone} centered on $x$ as the set of all vertices $y$ such that there is an oriented path from $y$ to $x$, but not to $\t x$.

\begin{figure}
	\centering
		\includegraphics[width=.50\textwidth]{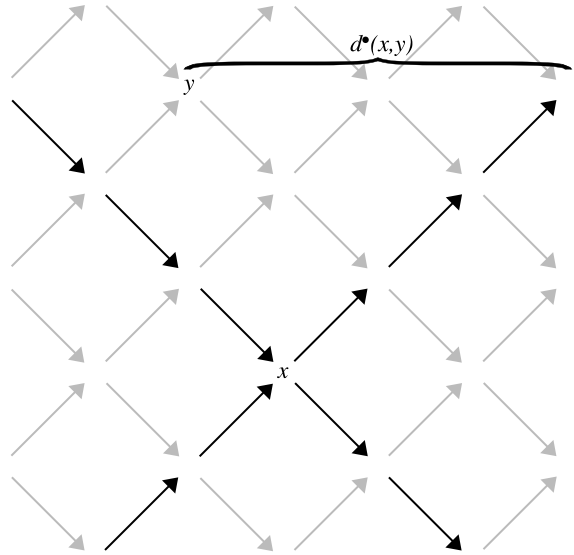}
	\caption{Light cones and light cone distance in $\bZ A_{\infty}^{\infty}$}
	\label{figure:LightCones}
\end{figure}

Let $y \in \bZ Q$ such that $\t^{-n} y$ lies on the right light cone centered on $x$, then we will say that the \emph{right light cone distance} $\r(x,y)$ is $n$.  Note that $\r(x,y)$ may be negative, and is not symmetric.  Fixing an $x$, the right light cone distance $\r(x,y)$ determines which object we take from the $\t$-orbit of $y$ (see for example Figure \ref{figure:LightCones}).

In Proposition \ref{proposition:rGeq0} we show that in order for a full subquiver $Q'$ of $\bZ Q$ to be a section, it suffices that $Q'$ meets every $\t$-orbit of $\bZ Q$ at least once, and that for every two vertices $x,y \in Q'$ both $\r(x,y)$ and $\r(y,x)$ are positive.  Graphically, these last conditions mean that $y$ lies ``in between'' the left and right light cones centered on $x$ (as is for example the case in Figure \ref{figure:LightCones}).

Another useful property of the right light cone distance is that one may see whether a certain section is strongly locally finite or not (Proposition \ref{proposition:PathFinite}).

Thus for the quiver $Q$, we pick any vertex $x \in \bZ Q$ and consider the left and right light cones centered on $x$.  In every $\t$-orbit, we choose a vertex ``in the middle'' between the left and right light cone centered on $x$ (as is illustrated in Figure \ref{figure:TiltedQuiver}).  Using properties of $\r$ we may then show that the constructed subquiver of $\bZ Q$ is a strongly locally finite section, completing the proof of Theorem \ref{theorem:Introduction}.

\textbf{Acknowledgment}  The authors would like to thank Idun Reiten and Sverre Smal\o\ for many useful discussions and helpful ideas.  We thank Michel Van den Bergh for his comments on an earlier version of the paper.  The second author also gratefully acknowledges the hospitality and support of the Max-Planck-Institut f\"ur Mathematik in Bonn and the Norwegian University of Science and Technology.

%% file: Preliminaries.tex
\section{Preliminaries}\label{section:preliminaries}

\subsection{Quivers} A quiver $Q$ is a 2-tuple $(Q_0,Q_1)$ of sets where the elements of $Q_0$ are the vertices, and $Q_1$ consists of arrows between those vertices.  We will often write $x \in Q$ for a vertex $x$, when we mean $x \in Q_0$.

An \emph{(oriented) path} between two vertices $x,y \in Q$ is a sequence $x=x_0, x_1, \ldots, x_{n-1}, x_n = y$ such that there is an arrow $x_i \to x_{i+1}$, for all $i \in \{0, 1, \ldots n-1\}$.  An (oriented) cycle is a nontrivial path from a vertex to itself.

We define \emph{unoriented paths} in an obvious way.  While we will often abbreviate ``oriented paths'' to ``paths'', in order to avoid confusion we will not abbreviate ``unoriented paths'' to ``paths''.

If there is an arrow $x \to y$ between two vertices, we say $x$ is a \emph{neighbor} of $y$ and vice versa.  If every vertex of $Q$ has only finitely many neighbors, we say $Q$ is \emph{locally finite}.  If $Q$ does not contain a subquiver of the form $\cdot \to \cdot \to \cdots$ or $\cdots \to \cdot \to \cdot$ (called \emph{rays} and \emph{corays}, respectively), we will say $Q$ is \emph{path finite}.

A connected, locally finite and path finite quiver has been called \emph{strongly locally finite} in \cite{ReVdB02}.  Hence a quiver $Q$ is strongly locally finite if and only if all indecomposable projectives and injectives representations have finite length.

\subsection{Stable translation quivers}  A \emph{stable translation quiver} is a quiver $T = (T_0,T_1)$ together with a bijection $\t : T_0 \longrightarrow T_0$, such that for all vertices $x,y \in T_0$ the number of arrows from $y$ to $x$ is equal to the number of arrows from $\tau x$ to $y$.

With a quiver $Q$ we will associate a stable translation quiver $\mathbb{Z} Q$ as the quiver with vertices $\mathbb{Z} Q_0 = \{ (n,x) \mid n \in \mathbb{Z}, x \in Q_0 \}$ and arrows as follows: the number of arrows $(i,x) \to (j,y)$ is equal to the number of arrows $x \to y$ in $Q$ if $i=j$, to the number of arrows $y \to x$ if $j = i+1$, and zero otherwise.  The translation is given by $\t: (i,x) \mapsto (i-1,x)$.

It is easy to see that $\bZ Q$ is locally finite or contains no oriented cycles if and only if $Q$ is locally finite or contains no oriented cycles, respectively.

A \emph{section} $Q$ of the stable translation quiver $\bZ Q$ is a connected full convex subquiver that meets each $\t$-orbit of $\bZ Q$ exactly once.  Thus $Q' \subseteq \bZ Q$ is a section if and only if the canonical injection may be lifted to an isomorphism $\bZ Q' \cong \bZ Q$ of stable translation quivers.

An equivalent formulation is given by \cite{ReVdB02}: a subquiver $Q'$ of the stable translation quiver $\bZ Q$ is a section if and only if $Q'$ meets every $\t$-orbit of $\bZ Q$ exactly once, and if $x \in Q'$ and $x\to z$ is an arrow in $\bZ Q$ then either $z\in Q'$ or $\tau z \in Q'$, and when $z\to x$ is an arrow in $\bZ Q$ then either $z\in Q'$ or $\tm z \in Q'$.

A \emph{sectional path} in a stable translation quiver is an oriented path $A_0 \to A_1 \to \cdots \to A_n$ such that $A_i \not= \t A_{i+2}$, for all $i \in \{0, \ldots, n-2\}$.

We will mostly be interested in stable translation quivers of the form $\bZ Q$ where there will only be a finite number of paths between two vertices. The following proposition gives some equivalent formulations.

\begin{proposition}\label{prop:SecFinite}
Let $Q$ be a connected quiver.  The following statements are
equivalent.
\begin{enumerate}
\item
The quiver $Q$ is locally finite and there are only finitely many
sectional paths between any two vertices of $\mathbb{Z} Q$.
\item
There are only finitely many (possibly non-sectional) paths between
any two vertices in $\mathbb{Z} Q$.
\item
For every vertex $x \in \mathbb{Z} Q$ there are only finitely many
paths from $x$ to $\tau^{-n} x$ in $\mathbb{Z} Q$ for all $n \in
\mathbb{N}$.
\item
There is a vertex $x \in \mathbb{Z} Q$ such that there are only
finitely many paths from $x$ to $\tau^{-n} x$ in $\mathbb{Z} Q$ for
all $n \in \mathbb{N}$.
\end{enumerate}
\end{proposition}
\begin{proof}
\begin{enumerate}
\item[$(1 \Rightarrow 2)$]
Seeking a contradiction to the assumptions in (1), we will assume we
may choose $x$ and $y$ such that there are infinitely many paths
from $x$ to $y$. Without loss of generality, we may assume $x$
has coordinates $(0,v_{x})$ and $y$ has coordinates $(n,v_{y})$,
where $v_x$ and $v_y$ are vertices in $Q$ and $n \geq 0$.

Since there are finitely many sectional paths from $x$ to
$y$, an infinite number of the paths between $x$ and $y$ must be
non-sectional. If $x\not=\tau y$ then we may turn a non-sectional
path into a non-trivial path from $x$ to $\tau y$ by replacing a part $A_{i-2} \to \t A_{i+1} \to A_i \to A_{i+1} \to A_{i+2}$ by $A_{i-2} \to \t A_{i+1} \to \t A_{i+2}$.

Since the paths from $x$ to $y$ have finite length and $Q$ is locally finite, only finitely many
different paths will be turned into the same one by this procedure, thus there are infinitely many paths from
$x$ to $\tau y$. Repeating this process shows that we either have
infinitely many paths from $x$ to $\tau^{n+1} y$ or infinitely many
paths from $x$ to $\tau^- x$.

The coordinates of $\tau^{n+1} y$ are $(-1,v_{y})$, and as such there may be
no paths from $x = (0,v_{x})$ to $\tau^{n+1} y$.

Therefore assume there are infinitely many paths from $x$ to $\tau^-
x$.  Since $Q$ is locally finite, there may only be a finite number
of paths from $x$ to $\tau^{-}x$ of length 2. 

All paths from $x$ to $\tau^{-}x$ not of length 2 are sectional, since otherwise we may
turn them into paths from $x$ to $x$ by replacing a part $A_{i-2} \to \t A_{i+1} \to A_i \to A_{i+1} \to A_{i+2}$ by $A_{i-2} \to \t A_{i+1} \to \t A_{i+2}$, as before.  Such a path from $x$ to $x$ is necessarily sectional.  By concatenating this cycle with itself, we obtain an infinite number of sectional paths from $x$ to $x$, a contradiction.

Hence we know there are infinitely many sectional paths from $x$ to $\tau^- x$, a contradiction to the assumption in (1).

\item[$(2 \Rightarrow 1)$]
  There is a finite number of paths between $x$ and $\tm x$ such that $Q$ is locally finite.  The claim about sectional paths is trivial.

\item[$(2 \Rightarrow 3 )$] Trivial.

\item[$(3 \Rightarrow 4 )$] Trivial.

\item[$(4 \Rightarrow 2)$]
Seeking a contradiction, assume there are infinitely many paths from
a vertex $y$ to a vertex $z$ of $\mathbb{Z} Q$.  Since $Q$ is
connected, there is a path from $x$ to $\tau^n y$ for an $n \in
\mathbb{Z}$. For the same reason there is a path from $\tau^n z$ to
$\tau^m x$ for an $m \in \mathbb{Z}$.  Composition gives a path from
$x$ to $\tau^m x$, hence $m \in -\mathbb{N}$.  Since there are
infinitely many paths from $y$ to $z$, composition gives infinitely
many paths from $x$ to $\tau^{m} x$, a contradiction to the
assumption in (4).
\end{enumerate}
\end{proof}

%% file: Distances.tex
\section{Light Cone and Round Trip Distance}\label{section:Distances}

In this section, we will introduce some tools that will help us to find and discuss sections in stable translation quivers of the form $\bZ Q$.

\subsection{Right light cone distances}

Let $Q$ be a quiver.  In $\bZ Q$ we define the \emph{(right) light cone}\index{right light cone!in a stable translation quiver} centered on a vertex $x \in \bZ Q$ as the set of all vertices $y$ such that there is a path from $x$ to $y$ but not to $\t y$.  It is clear that the right light cone intersects a $\t$-orbit in at most one vertex.  If $Q$ (and hence $\bZ Q$) is connected, then the right light cone intersects each $\t$-orbit in exactly one vertex.
\index{light cone|see{right light cone}}
\begin{figure}
	\centering
	\psfrag{x}[][]{$x$}
		\includegraphics[width=0.50\textwidth]{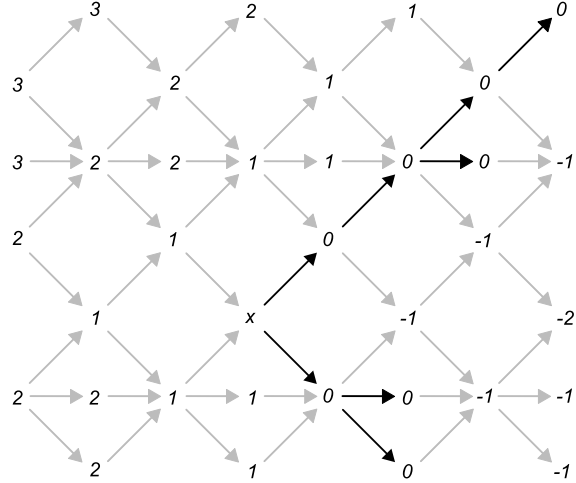}
	\caption{A stable translation quiver with the (right) light cone centered on $x$ and the corresponding right light cone distances}
	\label{fig:RightNumbered}
\end{figure}

Let $y \in \bZ Q$.  If $\t^{-n} y$ lies on the right light cone centered on $x$, then we will say that the \emph{right light cone distance}\index{right light cone distance!for stable translation quivers} $\r(x,y)$ is $n$.  If no such $n$ exists, we define $d^\bullet (x,y)=\infty$.  If $Q$ (and hence $\bZ Q$) is connected then the right light cone distance $d^\bullet (x,y)$ is finite for all vertices $x, y \in \mathbb{Z} Q$.

The following lemma is obvious.

\begin{lemma}\label{lemma:LightConeDistance}
For all $x,y \in \bZ Q$, we have $\r(x,\t^n y) = \r(x,y) + n$.
\end{lemma}

Note that $\r(x,y)$ may be negative, and that the function $\r$ is not symmetric.  The following lemma shows the right light cone distance satisfies the triangle inequality.

\begin{lemma}\label{lemma:rTriangle}
For all vertices $x,y,z \in \mathbb{Z} Q$ we have
$$d^\bullet (x,z) \leq d^\bullet (x,y) + d^\bullet (y,z)$$
\end{lemma}
\begin{proof}
Assume $d^\bullet (x,y) = n$ and $d^\bullet (y,z) = m$, thus there
are paths from $x$ to $\tau^{-n} y$ and from $\tau^{-n}y$ to
$\tau^{-n-m} z$. Composition gives a path from $x$ to $\tau^{-n-m}
z$, hence $d^\bullet (x,z) \leq n+m$.  If either $d^\bullet (x,y)$
or $d^\bullet (y,z)$ is infinite, then the inequality is trivial.
\end{proof}

There is a natural embedding $\epsilon \colon Q \hookrightarrow
\mathbb{Z} Q$ induced by the map $\epsilon (x) = (0,x)$. Let $x$ and
$y$ be vertices of $Q$, then we define the \emph{right light cone
distance}\index{quiver!right light cone distance}\index{right
light cone!distance!in a quiver} $d^\bullet_Q(x,y)$ between $x$ and
$y$ as the distance $d^\bullet((0,x),(0,y))$.

An equivalent way to describe $d^\bullet_Q (x,y)$ intrinsically on $Q$ is as the minimal number of arrows traversed in the opposite direction over all unoriented paths from $x$ to $y$.

\begin{proposition}\label{proposition:Hemimetric}
Let $Q$ be a connected quiver, then $\rQ$ defines a hemimetric on $Q$, i.e. for all $x,y,z \in Q$ we have
\begin{enumerate}
  \item $\rQ(x,y) \geq 0$,
  \item $\rQ(x,x) = 0$,
  \item $\rQ(x,z) \geq \rQ(x,y) + \rQ(y,z)$.
\end{enumerate}
If furthermore $Q$ does not have oriented cycles, then we may strengthen $(2)$ to
\begin{enumerate}
  \item[(2')] $\mbox{($\rQ(x,y) = 0$ and $\rQ(y,x) = 0$)}\Longleftrightarrow x=y$.
\end{enumerate}
\end{proposition}

\begin{proof}
This follows directly from the definition of $\rQ$ and Lemma \ref{lemma:rTriangle}.
\end{proof}

\begin{proposition}\label{proposition:r=0_-1}
If $x \to y$ is an arrow in $\mathbb{Z} Q$ for a quiver $Q$, then $d^\bullet (x,y)=0$ or
$d^\bullet (x,y)=-1$.  Furthermore $Q$ has no oriented cycles if and only $d^\bullet (x,y)=0$ for all arrows $x \to y$.
\end{proposition}

\begin{proof}
By the definition of $d^\bullet (x,y)$, and since there is a path from $x$ to $y$, we have $d^\bullet (x,y) \leq 0$.

From the arrow $x \to y$ we easily obtain an arrow $\tau^2 y \to \tau x$.  A path $x \to \tau^n y$ for $n \geq 2$
would produce a path from $x$ to $\tau x$ by concatenation with a path from $\tau^n y$ to $\tau^2 y$ and the arrow from $\tau^2 y$ to $\tau x$.  From the definition of $\mathbb{Z} Q$ we see such a path does not occur, hence $d^\bullet (x,y) \geq -1$.  This shows $d^\bullet (x,y)=0$ or $d^\bullet (x,y)=-1$.

If $\d^\bullet(x,y) = -1$, then there is a path from $x$ to $\t y$.  The arrow $x \to y$ yields an arrow $\t y \to x$ and we obtain a cycle in $\bZ Q$.  This implies there is a cycle in $Q$ as well.

Finally, assume $Q$ admits a cycle, and let $x \to y$ be an arrow occurring in this cycle.  This implies there is also an arrow $y \to \tm x$.  Since there is a path from $y$ to $x$, we know $\r(y,x) \leq 0$, and hence $\r(y, \tm x) \leq -1$ in $\bZ Q$.
\end{proof}

\begin{example}
Let $Q = \tilde{A}_1$ with cyclic orientation, and let $x,y \in \bZ Q$ as follows
$$\xymatrix@R=5pt{ & \bullet \ar@/^/[dd] \ar '[rd] [rrdd] && \bullet \ar@/^/[dd] \ar[rrdd] && y \ar@/^/[dd] \ar '[rd] [rrdd] && \bullet  \ar@/^/[dd] &\\
\cdots & & & & & && & \cdots \\
&\bullet \ar@/^/[uu] \ar[rruu] && x \ar@/^/[uu] \ar '[ru] [rruu] && \bullet \ar@/^/[uu] \ar[rruu] && \bullet  \ar@/^/[uu] &
}$$
Then $\r(x,y)=-1$.
\end{example}

\begin{figure}
	\centering
		\psfrag{x}[][]{$x$}
		\includegraphics[width=0.50\textwidth]{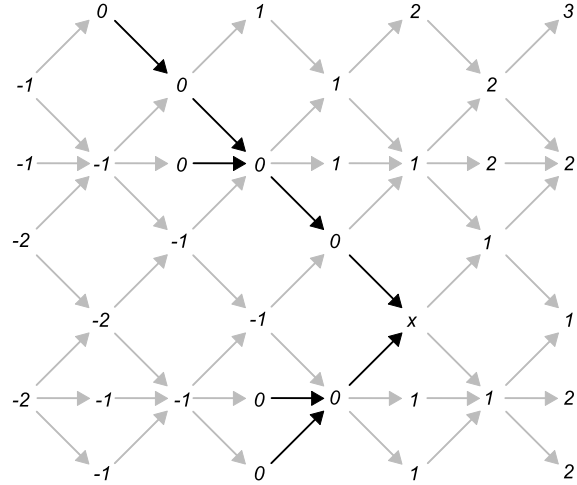}
		\caption{A stable translation quiver with the left light cone centered on $x$ and the corresponding left light cone distances}
	\label{fig:LeftNumbered}
\end{figure}

In addition to the right light cone distance one may also define a left light cone and a left light cone distance\index{left
light cone!distance} $d_\bullet \colon \mathbb{Z} Q \times \mathbb{Z}
Q \to \mathbb{Z} \cup \{ \infty \}$ dually (see Figure \ref{fig:LeftNumbered}), but since $\l (x,y) =
\r (y,x)$, the left light cone distance is essentially superfluous.

\subsection{Round Trip Distances}
\begin{figure}
	\centering
		\includegraphics[width=0.50\textwidth]{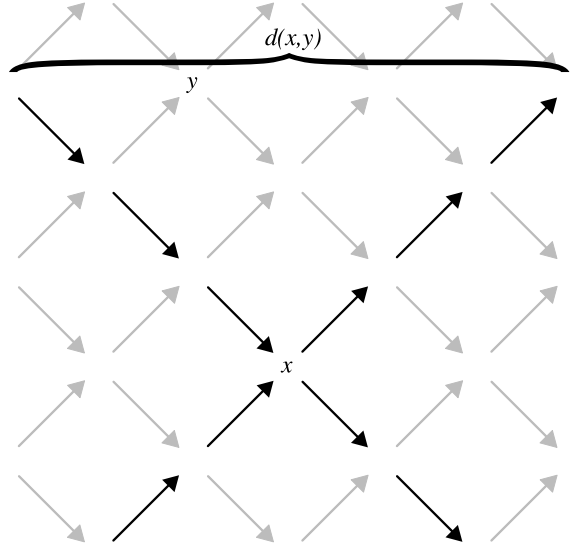}
	\caption{Light cones and round trip distance in $\bZ A^\infty_\infty$}
	\label{fig:RoundTrip}
\end{figure}

For two vertices $x, y \in \mathbb{Z} Q$, we define the \emph{round
trip distance}\index{round trip distance!in a translation
quiver}\index{translation!quiver!round trip distance} $d(x,y)$ as
$$d(x,y) = d^\bullet (x,y) + d^\bullet (y,x).$$
It is an immediate consequence of the definition that $d (x,y)$ is
the least integer $n$ such that there is a path in $\mathbb{Z} Q$
from $x$ to $\tau^{-n} x$ that contains exactly one vertex from the
$\tau$-orbit of $y$, namely $\tau^{-d^\bullet (x,y)} y$.

Let $x$ and $y$ be vertices of $Q$, then we define the \emph{round
trip distance}\index{round trip distance!in a
quiver}\index{quiver!round trip distance} $d_Q(x,y)$ between $x$ and
$y$ as the distance $d((0,x),(0,y))$ where $(0,x)$ and $(0,y)$ are
the vertices in $\mathbb{Z} Q$ corresponding to $x$ and $y$ under
the natural embedding $Q \hookrightarrow \mathbb{Z} Q$. Hence
\begin{align*} d_Q(x,y) &= d((0,x),(0,y)) \\ &= d^\bullet((0,x),(0,y)) +
d^\bullet((0,y),(0,x)) = d_Q^\bullet (x,y) + d_Q^\bullet (y,x)
\end{align*} 

As with $d^\bullet_Q$, we may describe $d_Q(x,y)$ intrinsically. If $x$ and $y$ are vertices of $Q$, then $d_Q(x,y)$
is the least number of arrows traversed in the opposite direction on
a path from $x$ to itself passing through $y$.

The next proposition shows that the round trip distance $\dQ$ defines a distance on the vertices of $Q$ when $Q$ is without oriented cycles.  If $Q$ has oriented cycles, then $\d$ merely defines a pseudodistance (i.e. satisfies conditions $(1)$ to $(4)$ below).

\begin{figure}
	\centering
	\psfrag{x}[][]{$x$}
		\includegraphics[width=0.58\textwidth]{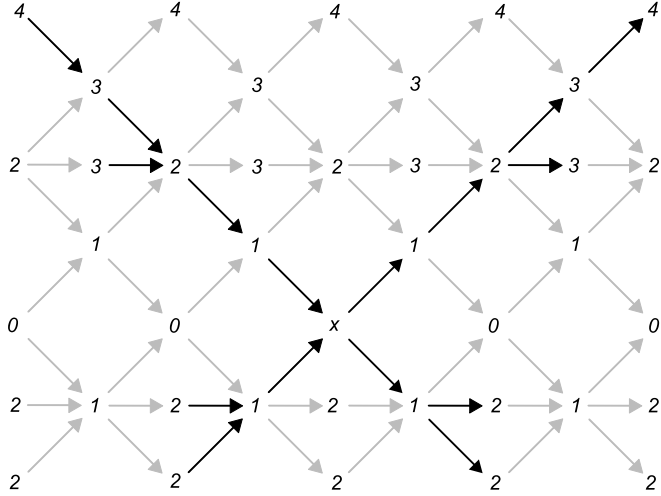}
	\caption{An example of a stable translation quiver where the left and right light cone centered on $x$ have been marked.  Every vertex is labeled with the round trip distance from $x$.}
	\label{fig:RoundTripNumbered}
\end{figure}

\begin{proposition}\label{proposition:d_is_distance}
Let $Q$ be a connected quiver, then for all $x,y,z \in Q$ we have
\begin{enumerate}
\item $d_Q(x,y) \geq 0$
\item $d_Q(x,x)=0$
\item $d_Q(x,y) = d_Q(y,x)$
\item $d_Q(x,z) \leq d_Q(x,y) + d_Q(y,z)$
\end{enumerate}
Furthermore, if $Q$ has no oriented cycles then we may strengthen $(2)$ to
\begin{itemize} \item[(2')] $d_Q(x,y) = 0 \Leftrightarrow x=y$
\end{itemize}

\end{proposition}

\begin{proof}
The first three properties follow directly from the definition of $\dQ$, while the triangle inequality follows from Lemma \ref{lemma:rTriangle}.  Furthermore, if $\dQ(x,y)=0$, then $x$ and $y$ lie on the same oriented cycle in $Q$.  This proves the last assertion.
\end{proof}

\subsection{Round Trip Distance Spheres for Quivers}
For a vertex $x$ in a quiver $Q$ we define the \emph{round trip
distance spheres}\index{round trip distance spheres} $S_Q(x,n)$
where $n \in \mathbb{N}$, as the sets 
$$S_Q(x,n) = \{ y \in Q \mid d_Q(x,y) = n \}.$$
Similarly we define the \emph{right light cone
sphere}\index{right light cone!sphere} and the \emph{left light cone
sphere}\index{left light cone!sphere} as 
$$\mbox{$S^\bullet_Q(x,n) = \{ y \in Q \mid d^\bullet_Q(x,y) = n \}$ and $S_\bullet^Q(x,n) = \{ y \in Q \mid d^\bullet_Q(y,x) = n \}$}$$ respectively.

We may now extend Proposition \ref{prop:SecFinite}.

\begin{proposition}\label{prop:dnFinite}
Let $Q$ be a connected quiver.  The following statements are
equivalent.
\begin{enumerate}
\item
The quiver $Q$ is locally finite and there are only finitely many
sectional paths between any two vertices of $\mathbb{Z} Q$.
\item
There are only finitely many (possibly non-sectional) paths between
any two vertices in $\mathbb{Z} Q$.
\item
For every vertex $x \in \mathbb{Z} Q$ there are only finitely many
paths from $x$ to $\tau^{-n} x$ in $\mathbb{Z} Q$ for all $n \in
\mathbb{N}$.
\item
There is a vertex $x \in \mathbb{Z} Q$ such that there are only
finitely many paths from $x$ to $\tau^{-n} x$ in $\mathbb{Z} Q$ for
all $n \in \mathbb{N}$.
\item
The quiver $Q$ is without oriented cycles, and for all $x \in Q$ and
$n \in \mathbb{N}$ the round trip distance sphere $S_Q(x,n)$ is
finite.
\item
The quiver $Q$ is without oriented cycles, and there is an $x \in Q$ such that the round trip distance sphere $S_Q(x,n)$ is finite for all $n \in \bN$.
\end{enumerate}
\end{proposition}
\begin{proof}
\begin{enumerate}

\item[$(3 \Rightarrow 5)$]
Since an oriented cycle involving $x$ would give infinitely many paths from $x$ to $x$, we know $Q$ is without oriented cycles.

Since every vertex $y \in S_Q(x,n)$ has a $\t$-shift in $\bZ Q$ lying on a path from $x$ to $\t^{-n} x$, and there are only finitely many such paths, it is clear $S_Q(x,n)$ must be finite.

\item[$(5 \Rightarrow 2)$]
For every $y$ on a path from $x$ to $\t^{-n} x$, we know $\d(x,y) \leq n$.  Since $S_Q(x,i)$ is finite for all $i \leq n$, there may only be finitely many paths from $x$ to $\t^{-n} x$.

\item[$(5 \Leftrightarrow 6)$] This follows directly from the triangle inequality.
\end{enumerate}
\end{proof}

%% file: Main.tex
\section{Existence of Strongly Locally Finite Sections}

We will now turn our attention to finding strongly locally finite sections in translation quivers of the form $\bZ Q$.  To do this we will use the right light cone distance and the round trip distance introduced in Section \ref{section:Distances}.

First, we will give a characterization of strongly locally finite quivers using the right and left light cone distances.

\begin{proposition}\label{proposition:PathFinite}
Let $Q$ be a connected quiver.  Then $Q$ is strongly locally finite if and only if $Q$ has no oriented cycles and for any $x \in Q$ all spheres $S^{\bullet}_Q(x,n)$ and $S^Q_{\bullet}(x,n)$ are finite for all $n \in \mathbb{N}$.
\end{proposition}

\begin{proof}
First, assume $Q$ is strongly locally finite.  Since $Q$ is then path finite, it is clear that $Q$ does not have oriented cycles.  Seeking a contradiction, we will assume there to be an $m \in \bN$ such that $S^{\bullet}_Q(x,m)$ is infinite for a certain vertex $x \in Q$.  Let $m$ be the smallest such integer; since $Q$ is path finite, we know $m \geq 1$.

For every $y \in S^{\bullet}_Q(x,m)$ we fix an unoriented path from $x$ to $y$ with exactly $m$ arrows in the opposite direction.  Following such an unoriented path from $x$ to $y$, the right light cone distance will be increasing.  Let $z$ be the first vertex encountered on this unoriented path with $\r(x,z) = m$.

Such a vertex $z$ admits an oriented path to $y$ and is a neighbor of a vertex in $S^{\bullet}_Q(x,m-1)$.  Since this last set is finite and $Q$ is locally finite, it is clear that there are only finitely many vertices $z$.  Hence one of these vertices admits oriented paths to an infinite number of vertices in $S^{\bullet}_Q(x,m)$.  Since $Q$ is locally finite, we conclude that $Q$ has rays.  A contradiction.

Dually, one shows $S_{\bullet}^Q(x,n)$ is finite for all $n \in \mathbb{N}$.

For the other implication, assume $Q$ has no oriented cycles and for a certain $x \in Q$ all spheres $S^{\bullet}_Q(x,n)$ and $S^Q_{\bullet}(x,n)$ are finite for all $n \in \mathbb{N}$.

Let $y \in Q$ be any vertex.  For all neighbors $z \in Q$ of $y$, we have either $\rQ(y,z)=0$ if there is an arrow $y \to z$ or $\rQ(y,z)=1$ if there is an arrow $z \to y$.  Using the triangle inequality, we find 
$$\rQ(x,z) \leq \rQ(x,y) + \rQ(y,z) \leq \rQ(x,y) + 1.$$
Since $S^{\bullet}_Q(x,n)$ is finite for all $n \in \bN$, we see that $y$ may only have a finite number of neighbors, hence $Q$ is locally finite.

We will now proceed by proving that $Q$ is path finite.  Assume $Q$ has a ray $a_0 \to a_1 \to \cdots$ as subquiver.  For $i \geq 0$, the triangle inequality gives
$$d^\bullet_Q(x, a_{i+1}) \leq d^\bullet_Q(x,a_i) + d^\bullet_Q(a_i, a_{i+1}) = d^\bullet_Q(x, a_i)$$
since $d^\bullet_Q(a_i, a_{i+1}) = 0$, hence the sequence $(d^\bullet_Q(x, a_i))_{i \in \mathbb{N}}$ must stabilize, giving an infinite set $S^{\bullet}_Q(x,m)$ for an $m \leq d^\bullet_Q(x, a_0)$.  Thus $Q$ may not have a ray as a subquiver.

Dually, one finds that $Q$ may not have a coray as subquiver.
\end{proof}

The next result gives necessary and sufficient conditions for $Q'$ to be a section of $\mathbb{Z} Q$ using the right light cone distance.

\begin{proposition}\label{proposition:rGeq0}
Let $Q'$ be a full subquiver of the stable translation quiver $\mathbb{Z} Q$ that meets every $\tau$-orbit exactly once. Then $Q'$ is a section if and only if $d^\bullet(x,y) \geq 0$ for all vertices $x,y\in Q'$.
\end{proposition}

\begin{proof}
We will first check that, if $d^\bullet(x,y) \geq 0$ for all vertices $x,y\in Q'$, then $Q'$ is a section.  We need to show that for every arrow
$x\to z$ in $\mathbb{Z} Q$ with $x \in Q'$ either $z\in Q'$ or $\tau
z \in Q'$, and for every arrow $z\to x$ in $\mathbb{Z} Q$ with $x
\in Q'$ either $z\in Q'$ or $\tm z \in Q'$. We will only show the
first part, the second is similar.

So let $x \in Q'$.  Since there is an arrow $x\to z$ in $\bZ Q$, we know $\r(x,z) \leq 0$, thus the object of the $\t$-orbit of $z$ belonging to $Q'$ has to be of the form $\t^n z$ with $n \geq 0$.

An arrow $x\to z$ induces an arrow $\t z \to x$, hence $\r(\t z, x) \leq 0$ and thus $n \leq 1$.  We conclude that either $z$ or $\t z$ belongs to $Q'$.

Conversely, let $Q'$ be a section of $\bZ Q$ and let $x,y \in Q'$.  Since the injection $Q' \subseteq \bZ Q$ lifts to an isomorphism $\bZ Q' \to \bZ Q$ of translation quivers, Proposition \ref{proposition:Hemimetric} yields $d^\bullet (x,y) = d^\bullet_{Q'}(x,y) \geq 0$.
\end{proof}

\begin{example}
Let $x$ be a vertex of the stable translation quiver $\mathbb{Z} Q$.  Using triangle inequality, one easily verifies that the right light cone $\Sr x 0$ and the left light cone $\Sl x 0$ are both sections of $\bZ Q$.
\end{example}

We now come to the main result of this section.

\begin{theorem}\label{theorem:Main}
Let $Q$ be a connected quiver.  The following statements are equivalent.
\begin{enumerate}
\item
The quiver $Q$ is locally finite and there are only finitely many sectional paths between any two vertices of $\mathbb{Z} Q$.
\item
There are only finitely many (possibly non-sectional) paths between any two vertices in $\mathbb{Z} Q$.
\item
For every vertex $x \in \mathbb{Z} Q$ there are only finitely many paths from $x$ to $\tau^{-n} x$ in $\mathbb{Z} Q$ for all $n \in \mathbb{N}$.
\item
There is a vertex $x \in \mathbb{Z} Q$ such that there are only finitely many paths from $x$ to $\tau^{-n} x$ in $\mathbb{Z} Q$ for all $n \in \mathbb{N}$.
\item
The quiver $Q$ is without oriented cycles, and for all $x \in Q$ and $n \in \mathbb{N}$ the round trip distance sphere $S_Q(x,n)$ is finite.
\item
The quiver $Q$ is without oriented cycles, and there is an $x \in Q$ such that the round trip distance sphere $S_Q(x,n)$ is finite, for all $n \in \bN$.
\item
The translation quiver $\mathbb{Z} Q$ has a strongly locally finite section.
\end{enumerate}
\end{theorem}

\begin{proof}
The first 6 points are equivalent by Proposition \ref{prop:dnFinite}.
\begin{itemize}
\item[$(5 \Rightarrow 7)$]
We will construct a section $Q'$ in $\mathbb{Z} Q$.  Start by fixing
a vertex $x$ in $\mathbb{Z} Q$.  From every $\tau$-orbit we will
choose a vertex $y$ to be in $Q'$ for which $d^\bullet(x,y) =
\left\lfloor \frac{d(x,y)}{2}\right\rfloor$, hence also $d^\bullet(y,x) = \left\lceil \frac{d(x,y)}{2}\right\rceil$, where $\lfloor \cdot \rfloor$ and $\lceil \cdot \rceil$ are the usual floor and ceiling functions, respectively.  We will use Proposition \ref{proposition:rGeq0} to show that the full subquiver $Q'$ picked in this way is a section of $\mathbb{Z} Q$.

Therefore we need to show that for all vertices  $y,z\in Q' \subset
\mathbb{Z} Q$, we have $d^\bullet_{Q'}(y,z) \geq 0$.  We will
consider two cases.  First, assume $d(x,z) - d(x,y) \geq 0$.
Starting with the triangle inequality, we have
\begin{eqnarray*}
d^\bullet(y,z) &\geq& d^\bullet(x,z) - d^\bullet(x,y) \\
 &=& \left\lfloor \frac{d(x,z)}{2}\right\rfloor - \left\lfloor \frac{d(x,y)}{2}\right\rfloor\\
&\geq& 0
\end{eqnarray*}
Next if $d(x,z) - d(x,y) \leq 0$, we have
\begin{eqnarray*}
d^\bullet(y,z) &\geq& d^\bullet(y,x) - d^\bullet(z,x) \\
&=& \left\lceil \frac{d(x,y)}{2}\right\rceil - \left\lceil \frac{d(x,z)}{2}\right\rceil\\
&\geq& 0
\end{eqnarray*}
Proposition \ref{proposition:rGeq0} then yields that $Q'$ is a
section of $\mathbb{Z} Q$.

To show that $Q'$ is path finite, we note that
$|S^\bullet_{Q'}(x,n)| = |S_Q(x,2n)| + |S_Q(x,2n+1)|$ and
$|S_\bullet^{Q'}(x,n)| = |S_Q(x,2n-1)| + |S_Q(x,2n)|$, so by
assumption the sets $S^\bullet_{Q'}(x,n)$ and $S_\bullet^{Q'}(x,n)$
are finite.  Since $Q$, and hence also $\bZ Q$, is locally finite and has no oriented cycles, we know that the same is true for $Q'$.  Proposition  \ref{proposition:PathFinite} now yields $Q'$ is path finite.

\item[$(7 \Rightarrow 5)$]
Let $Q'$ be a strongly locally finite section of $\bZ Q$.  We may assume there is a vertex $x \in \bZ Q$ lying in both $Q$ and $Q'$.  It is then clear that
$$|S_Q (x,n)| = |S_{Q'} (x,n)| = \left| \bigcup_{i+j = n} (S^\bullet_{Q'}(x,i) \cap S_\bullet^{Q'}(x,j)) \right|.$$
By Proposition \ref{proposition:PathFinite}, the right hand side is finite,
hence also the left hand side is finite.  Since $Q'$ is path finite,
it has no oriented cycles, so $Q$ is also without oriented cycles.
\end{itemize}
\end{proof}

\begin{example}

Let $Q$ be the quiver $A_\infty^\infty$ with linear orientation, thus 
$$Q: \cdots \to \cdot \to \cdot \to \cdot \to \cdot \to \cdot \to \cdot \to \cdots$$
It is easy to see that $Q$ satisfies statement $(6)$ in Theorem \ref{theorem:Main}.  After fixing a vertex $x$ of $\bZ Q$, the construction described in the proof of Theorem \ref{theorem:Main} gives a strongly locally finite quiver $Q'$ as in Figure \ref{figure:TiltedQuiver}, namely $Q'$ is an $A_\infty^\infty$-quiver with zig-zag orientation.
$$Q': \cdots \to \cdot \leftarrow \cdot \to \cdot \leftarrow \cdot \to \cdot \leftarrow \cdot \to \cdots$$

\begin{figure}
	\centering
	\psfrag{x}[][]{$x$}
		\includegraphics[width=0.58\textwidth]{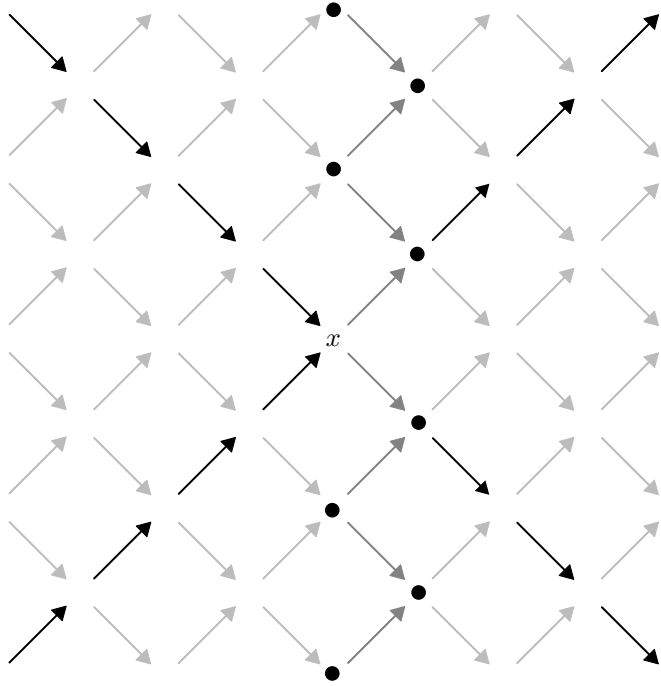}
	\caption{A stable translation quiver of the form $\bZ A_\infty^\infty$.  Here an object $x$ has been chosen and the associated left and right light cones are given by black arrows.  The vertices of the full subquiver $Q'$ constructed in the proof of Theorem \ref{theorem:Main} are indicated by `$\bullet$'. }
	\label{figure:TiltedQuiver}
\end{figure}
\end{example}

%% file: Categories.tex
\section{Application to Hereditary Categories with Serre Duality}

In this section, we apply Theorem \ref{theorem:Main} to the theory of $k$-linear abelian Ext-finite hereditary categories with Serre duality.  In this way, we contribute to an ongoing project to better understand these categories (cf. \cite{Happel01}, \cite{ReVdB02},  \cite{vanRoosmalen07}, \cite{vanRoosmalen06}).  Throughout, let $k$ be an algebraically closed field, and $\AA$ be a $k$-linear abelian Ext-finite hereditary category with Serre duality.  We start by recalling some definitions and a short discussion about sectional paths.

\subsection{Definitions}   Let $\mathcal{A}$ be an abelian $k$-linear category. We say $\mathcal{A}$ is \emph{hereditary} if $\Ext^2(X,Y) = 0$ for all $X,Y \in \mathcal{A}$ and is \emph{Ext-finite} if $\dim_k \Ext(X,Y) < \infty$ for all $X,Y \in \mathcal{A}$. 

We will say $\AA$ satisfies \emph{Serre duality} \cite{Bondal89} if there exists an auto-equivalence $F \colon D^b\mathcal{A} \to D^b\mathcal{A}$, called a \emph{Serre functor}, such that for all $X,Y \in \Db \AA$ there is an isomorphism
$$\Hom_{D^b \mathcal{A}}(X,Y) \cong \Hom_{D^b\mathcal{A}}(Y, F(A))^\ast$$
natural in $X$ and $Y$, where $(-)^\ast$ is the vector space dual.

It has been shown in \cite{ReVdB02} that $\AA$ has Serre duality if and only if $\Db \AA$ has Auslander-Reiten triangles; the Serre functor then coincides with $\t [1]$, where $\t$ is the Auslander-Reiten translation in $\Db \AA$. In particular, a hereditary category $\AA$ has Serre duality if and only if $\AA$ has Auslander-Reiten sequences and the Serre functor $F : \Db \AA \longrightarrow \Db \AA$ induces an equivalence between the category of projectives and the category of injectives of $\AA$.

The \emph{Auslander-Reiten quiver} of $\AA$ and $\Db \AA$ is defined as follows.  The set of vertices is $\ind \AA$ or $\ind \Db \AA$, respectively, and between two vertices $A, B$, there are $\dim_k \rad(A,B) / \rad^2(A,B)$ arrows from $A$ to $B$.  If $\AA$ is an abelian hereditary Ext-finite category with Serre duality, then the Auslander-Reiten quiver of $\Db \AA$ is a stable translation quiver with $\t = F[-1]$.

The full subquiver of the Auslander-Reiten quiver of $\AA$ spanned by all projective or injective objects in $\ind \AA$ is called the \emph{quiver of projectives} or \emph{injectives} of $\AA$, respectively.  A component of the Auslander-Reiten quiver of $\AA$ containing a projective object is called a \emph{preprojective component}.

If $\AA$ satisfies Serre duality, then the Auslander-Reiten component of $\Db \AA$ containing the projective quiver $Q$ is a stable translation quiver of the form $\bZ Q$ where the translation $\t$ is given by the Auslander-Reiten translation.  We will refer to this Auslander-Reiten component as the \emph{connecting component}.

Finally, we will say a component $\QQ$ of the Auslander-Reiten sequence of $\Db \AA$ is \emph{generalized standard} if $\infrad(X,Y)=0$ for all vertices $X,Y$ of $\QQ$.  In particular, if there is no oriented path from $X$ to $Y$ in the Auslander-Reiten quiver, then $\Hom(X,Y)=0$. 

\subsection{Sectional paths}  We will say a sequence $A_0 \to A_1 \to \cdots \to A_{n-1} \to A_n$ of irreducible maps between indecomposable objects in $\AA$ or $\Db \AA$ is \emph{sectional} if $A_{i} \not \cong \tau A_{i+2}$ for all $i \in \{0,\ldots,n-2\}$.  Note that a corresponding path in the Auslander-Reiten quiver is a sectional path.

\begin{proposition}\label{proposition:linearly_independent}
Let $\mathcal{A}$ be an abelian $\Ext$-finite category with Serre duality, then for every $X,Y \ind \Db \AA$ there may only be finitely many sectional paths from $X$ to $Y$.
\end{proposition}

\begin{proof}
Assume there are different sectional paths from $X$ to $Y$.  The arrows $A \to B$ in the Auslander-Reiten quiver of $\Db \AA$ give a basis of $\rad(A,B) / \rad^2(A,B)$.  With such a basis, we may associate linearly independent morphisms of $\rad(A,B)$.  Fix such a morphism for every arrow occurring in an above path from $X$ to $Y$ (if an arrow occurs more than once, we will associate the same morphism with it).

In this way, every sectional path corresponds to a morphism in $\Hom(X,Y)$.  We claim different sectional paths give rise to linearly independent morphisms.

Seeking a contradiction, consider the sectional sequences as depicted below
$$\xymatrix@R=0.5pc{ & A_0^0 \ar[r]^{f^0_1} & A_1^0 \ar[r] & \cdots \ar[r] & A_{n_0}^0 \ar[ddr]^{f^0_{n_0 + 1}} & \\ & A_0^1 \ar[r]_{f^1_1} & A_1^1 \ar[r] & \cdots \ar[r] & A_{n_1}^1 \ar[dr]_{f^1_{n_1 + 1}} & \\
X \ar[uur]^{f^0_0} \ar[ur]_{f^1_0} \ar[dr]_{f^m_0} & \vdots & \vdots && \vdots & Y \\ &
A^m_0 \ar[r]_{f^m_1} & A_1^m \ar[r] & \cdots \ar[r] & A_{n_m}^m \ar[ur]_{f^m_{n_m + 1}} &
}$$
such that there is a linear combination of the corresponding maps
$$\sum_{i=0}^m \alpha_i \left(\bigcirc_{k=0}^{n_i + 1} f^i_k\right)=0$$
where $\alpha_i \in k \setminus \{ 0 \}$, and where the correct order of composition is understood.  Keeping all paths that end with the morphism $f^0_{n_0 + 1}$ on the left hand side of the equation and moving the others to the right hand side, we find (possibly after renumbering the paths)
$$f^0_{n_0 + 1} \circ \left( \sum_{i=0}^{m_0} \alpha_i \left(\bigcirc_{k=0}^{n_i} f^i_k \right) \right) = - \sum_{i=m_0 + 1}^m \alpha_i \left(\bigcirc_{k=0}^{n_i + 1} f^i_k\right)$$

Denote $g_0 = \sum_{i=0}^{m_0} \alpha_i \left(\bigcirc_{k=0}^{n_i} f^i_k\right)$.  Considering the Auslander-Reiten triangle extending the irreducible maps $f_{n_i +1}^i: A_{n_i}^i \to Y$ gives following diagram.

$$\xymatrix{X \ar@{.>}[dd] \ar[dr]^{g_0} \ar[dddr] &&& \\
& A_{n_0}^0 \ar[dr]^{f^0_{n_0+1}} && \\ \tau Y \ar[ur] \ar[dr] && Y \ar[r] & \tau
Y [1] \\ & E_1 \ar[ur] && }$$

It follows that $g_0: X \to A^0_{n_0}$ factors through the map $\t Y \to A^0_{n_0}$. Likewise, we may split the compositions occurring in the definition of $g_0$ in two groups, with the group on the left hand side containing all the compositions ending in $f^0_{n_0}$.  After possibly renumbering the paths, we get
$$f^0_{n_0} \circ \left( \sum_{i=0}^{m_1} \alpha_i \left(\bigcirc_{k=0}^{n_i - 1} f^i_k \right) \right) = g_0 - \sum_{i=m_1 + 1}^{m_1} \alpha_i \left(\bigcirc_{k=0}^{n_i} f^i_k\right).$$

If we write $g_1 = \sum_{i=0}^{m_1} \alpha_i \left(\bigcirc_{k=0}^{n_i - 1} f^i_k \right)$, then we see from the following Auslander-Reiten triangle 
$$\xymatrix{X \ar@{.>}[dd] \ar[dr]^{g_1} \ar[dddr] &&& \\
& A_{n_0-1}^0 \ar[dr]^{f^0_{n_0}} && \\ \tau A_{n_0}^0 \ar[ur] \ar[dr] &&
A_{n_0}^0 \ar[r] & \tau A_{n_0}^0 [1] \\ & E_2 \ar[ur] && }$$
that $g_1$ factors through $\t A^0_{n_0} \to A^0_{n_0-1}$.

Since every considered path is different, iterating this procedure shows that the irreducible map $\alpha_0 f_0^0: X \to A_0^0$ factors through $E \to A_0^0$ as in the Auslander-Reiten triangle
$$\xymatrix{X \ar[dr]^{\alpha_0 id_X} \ar[dddr] &&& \\
& X \ar[dr]^{f^0_0} && \\ \tau  A_0^0 \ar[ur] \ar[dr] && A_0^0 \ar[r] & \tau
 A_1^0 [1] \\ & E \ar[ur] && }$$
which is clearly a contradiction.
\end{proof}

This implies that every stable component of the form $\bZ Q$ of the Auslander-Reiten quiver $\Db \AA$ satisfies the equivalent conditions of Theorem \ref{theorem:Main}.  In particular, we have the following corollary.

\begin{corollary}\label{corollary:ZQ}
Let $\AA$ be an abelian $\Ext$-finite $k$-linear category with Serre duality.  If a component of the Auslander-Reiten quiver of $\Db \AA$ is of the form $\bZ Q$, then $Q$ satisfies the equivalent conditions of Theorem \ref{theorem:Main}.
\end{corollary}

\subsection{Representations of strongly locally finite quivers}  Let $Q$ be a strongly locally finite quiver.  It is easy to see that this implies that there are only finitely many paths between two vertices of $Q$.

Let $\rep_k Q$ be the category of finitely presented $k$-representations of $Q$ and denote by $\PP$ and $\II$ the full subcategory of projectives and injectives, respectively.  With every vertex $x \in Q$ we may associate an indecomposable projective object $P_x$ and an indecomposable injective object $I_x$.  There is a canonical isomorphism $\nu_{x,y} : \Hom(P_x,P_y) \cong \Hom(I_x, I_y)$ since both vector spaces have the paths of from $y$ to $x$ as a basis.

We may consider the Nakayama functor $N : \PP \to \II$ where $N(P_x) = I_x$ and where the map $\Hom(P_x,P_y) \to \Hom(N(P_x), N(P_y))$ is given by the above isomorphism $\nu_{x,y}$.  The Nakayama functor is an equivalence of categories.

It follows from \cite[Lemma II.1.2]{ReVdB02} that the composition $$F:\Db \rep_k Q \cong \Kb \PP \stackrel{N}{\longrightarrow} \Kb \II \cong \Db \rep_k Q$$ is a right Serre functor.  Since $F$ is an equivalence, it is a Serre functor.  Hence $\rep_k Q$ satisfies Serre duality.

\subsection{Derived Equivalences} Assume that $\Db \AA$ is generated as a triangulated category by the connecting component $\CC$ and furthermore that the connecting component is generalized standard, thus $\infrad(X,Y)=0$ for all $X$ and $Y$ in $\ind \CC$.  If we denote the quiver of projectives in $\AA$ by $Q$, then the Auslander-Reiten quiver of $\CC$ will be a stable translation quiver of the form $\bZ Q$.

Let $X$ be any vertex in $Q$ and $X \to M_X \to \tm X \to X[1]$ be an Auslander-Reiten triangle.  Since $M_X$ has as many direct summands as $X$ has direct successors in $\bZ Q$, we see that $\bZ Q$ and hence also $Q$ must be locally finite.  Furthermore, it follows from Proposition \ref{proposition:linearly_independent} that there may be only finitely many sectional paths between any two vertices in $\bZ Q$, thus by Theorem \ref{theorem:Main} we know $\bZ Q$ admits a strongly locally finite section $Q'$.

Since $\CC$ is generalized standard and $\Ext(X,Y) \cong \Hom(Y,\t X)^*$, Proposition \ref{proposition:rGeq0} yields that the vertices of $Q'$ form a partial tilting set, i.e.\ $\Hom_{\Db \AA}(X,Y[n])=0$ for all $n \in \bZ \setminus \{0\}$ and all $X,Y \in Q'$.  It follows from \cite[Theorem 5.1]{vanRoosmalen06} that there is an exact fully faithful functor $i:\Db \rep_k (Q')^\circ \longrightarrow \Db \AA$ mapping $P_X$ to $X$, where $(Q')^\circ$ is the dual quiver of $Q'$.

Considering the exactness of $i$, and the connection between the Auslander-Reiten translation $\t$ and the Serre functor $F$, we may check that $i\circ F (P) \cong F\circ i (P)$ for all $P \in Q'$.  Hence the essential image of $i$ is closed under the action of the Serre functor of $\Db \AA$ and contains $\CC$.  Since $\CC$ generates $\Db \AA$, we conclude that $i$ is essentially surjective and thus an equivalence.

We have proven the following theorem.

\begin{theorem}\label{theorem:MainDerived}
Let $\AA$ be a $k$-linear abelian Ext-finite hereditary category with Serre duality.  Assume $\Db \AA$ is generated by its connecting component $\CC$ and that $\CC$ is generalized standard, then $\AA$ is derived equivalent to $\rep_k Q'$ where $Q'$ is strongly locally finite.
\end{theorem}

We now turn our attention to noetherian categories.  It has been shown in \cite[Theorem II.4.2]{ReVdB02} that in this case the category $\AA$ decomposes as a direct sum $\RR \oplus \QQ$ where $\RR$ has no nonzero preprojectives, nor nonzero preinjectives, and where $\QQ$ is generated by preprojectives.  Thus, when $\AA$ is a $k$-linear connected noetherian abelian Ext-finite hereditary category with Serre duality, saying that $\AA$ has at least one non-zero projective object is equivalent to saying that $\AA$ is generated by preprojectives.  In this case, $\Db \AA$ is generated by the connecting component.

We have following corollary as answer to the question posed in \cite{ReVdB02}.

\begin{corollary}\label{corollary:ReVdB02}
Let $\AA$ be a noetherian $k$-linear abelian Ext-finite hereditary category with Serre duality.  Assume $\AA$ has a non-zero projective object, then $\AA$ is derived equivalent to $\rep_k Q'$ where $Q'$ is strongly locally finite.
\end{corollary}

\begin{proof}
It has been shown in \cite[Proposition II.2.3]{ReVdB02} that the quiver of projectives $Q$ of $\AA$ is locally finite and does not contain any subquivers of the form $\cdot\to\cdot\to\cdot\to\cdots$

Since $\AA$ has a nonzero projective object, it is generated by preprojectives and hence $\Db \AA$ is generated by the connecting component.

We will show the connecting component $\CC$ is generalized standard.  Let $X,Y \in \ind \CC$ be with coordinates $(0,v_X)$ and $(n,v_Y)$, respectively, in the Auslander-Reiten quiver $\bZ Q$ of $\CC$.  Assume that $\infrad(X,Y) \not= 0$ and that $n$ has been chosen minimal with this property.

Consider the Auslander-Reiten triangle $Y \to M_Y \to \tm Y  \to Y[1]$.  There is at least one indecomposable summand of $Y_1$ of $M_Y$ such that $\infrad(X,Y_1) \not= 0$.  Due to the minimality of $n$, the coordinates of $Y_1$ in $\bZ Q$ must be $(n,v_{Y_1})$ where $v_{Y_1}$ is a direct successor of $v_Y$ in $Q$.  Iteration gives an infinite sequence $Y \to Y_1 \to Y_2 \to \cdots$ in $Q$, a contradiction.

We may now use Theorem \ref{theorem:MainDerived} to see that $\AA$ is derived equivalent to $\rep_k Q'$ where $Q'$ is strongly locally finite.
\end{proof}

  Our last result characterizes all quivers which can occur as quiver of projectives of an abelian hereditary category with Serre duality.
  
\begin{corollary}\label{corollary:Classification}
Let $Q$ be a quiver.  The quiver $Q$ satisfies the equivalent conditions of Theorem \ref{theorem:Main} if and only if there is an abelian hereditary category with Serre duality having $Q$ as its quiver of projectives.
\end{corollary}

\begin{proof}
One direction follows directly from Corollary \ref{corollary:ZQ}.  So let $Q$ be a quiver satisfying the equivalent conditions of Theorem \ref{theorem:Main} and let $Q'$ be a strongly locally finite section in $\bZ Q$.  Consider the hereditary abelian category $\AA = \rep_k (Q')^{\circ}$.  The category of projectives of $\AA$ is then given by $Q'$.

We may assume $Q$ is an infinite quiver, and in particular not Dynkin.  The required result then follows from \cite[Lemma II.3.4]{ReVdB02}.
\end{proof}

%% file: ProjectiveQuiver.bbl
\begin{thebibliography}{1}

\bibitem{Bondal89}
A.~I. Bondal and M.~M. Kapranov, \emph{Representable functors, {S}erre
  functors, and reconstructions}, Izv. Akad. Nauk SSSR Ser. Mat. \textbf{53}
  (1989), no.~6, 1183--1205, 1337.

\bibitem{Happel01}
Dieter Happel, \emph{A characterization of hereditary categories with tilting
  object}, Invent. Math. \textbf{144} (2001), no.~2, 381--398. \MR{MR1827736
  (2002a:18014)}
  
\bibitem{ReVdB02}
I.~Reiten and M.~Van~den Bergh, \emph{Noetherian hereditary abelian categories
  satisfying {S}erre duality}, J. Amer. Math. Soc. \textbf{15} (2002), no.~2,
  295--366 (electronic). \MR{MR1887637 (2003a:18011)}

\bibitem{Ringel02b}
Claus~Michael Ringel, \emph{A ray quiver construction of hereditary abelian
  categories with Serre duality}, Representations of algebra. Vol. II, BNU
  Press, 2002, pp.~396--416.
  
\bibitem{vanRoosmalen07}
Adam-Christiaan van Roosmalen, \emph{Abelian 1-{C}alabi-{Y}au categories}, Int.
  Math. Res. Not. IMRN (2008), no.~6, Art. ID rnn003, 20.

\bibitem{vanRoosmalen06}
\bysame, \emph{Classification of abelian hereditary directed categories
  satisfying {S}erre duality}, Trans. Amer. Math. Soc. \textbf{360} (2008),
  no.~5, 2467--2503.
  
\end{thebibliography}
